\documentclass[%
 reprint,
nofootinbib,
 amsmath,amssymb,
 aps,
]{revtex4-1}

\pdfoutput=1

\usepackage{graphicx}
\usepackage{dcolumn}
\usepackage{bm}

\usepackage{physics}
\usepackage{wrapfig}
\usepackage{hanging}
\usepackage{hhline}

\begin{document}

\widowpenalty10000
\clubpenalty10000

\title{Dodecatonic Cycles and Parsimonious Voice-Leading\\in the Mystic-Wozzeck Genus}
\thanks{I thank Professor Suzannah Clark for discussions during the preparation of this paper.}%

\author{Vaibhav Mohanty}
\email{E-mail: mohanty@college.harvard.edu}
\affiliation{%
 Quincy House, Harvard University, Cambridge, MA 02138
}%

\count\footins = 10000


\begin{abstract}
This paper develops a unified voice-leading model for the genus of mystic and Wozzeck chords. These voice-leading regions are constructed by perturbing symmetric partitions of the octave, and new Neo-Riemannian transformations between nearly symmetric hexachords are defined. The behaviors of these transformations are shown within visual representations of the voice-leading regions for the mystic-Wozzeck genus.
\end{abstract}

\maketitle

\section*{I. Introduction}
\indent [1.1] In a footnote in his 1996 article, Richard Cohn mentions that it should be possible to understand voice-leading parsimony of chords of the pitch-class sets 4-27 and 6-34 just as he has shown for 3-11 (the consonant triads). Adrian Childs' 1998 paper full describes the theory for 4-27 in his paper two years later. In \textit{Audacious Euphony: Chromatic Harmony and the Triad's Second Nature}, Cohn (2012) reviews both the $n = 3$ and $n = 4$ cases, and he briefly discusses voice-leading parsimony in the mystic and Wozzeck chords for the $n = 6$ case, but he does not fully develop a unified voice-leading model. The central focus of this paper is to formalize the $n = 6$ voice-leading regions, developing a ``centipede'' region analogous to the Weitzmann waterbug and Boretz spider as well as a ``dodecatonic'' region---in analogy with the hexatonic and octatonic regions---for the $T_n/T_nI$ set class of nearly symmetric hexachords.

[1.2] In \textit{Audacious Euphony}, Richard Cohn (2012) unifies hexatonic cycles and Weitzmann regions/waterbugs into a single model of voice-leading for all 24 major and minor triads (\textit{i.e.}, nearly symmetric chords of cardinality $n = 3$). By analogy, he also constructs a unified geometric representation for dominant seventh and fully diminished seventh chords (the $n = 4$ case), combining Childs' (1998) octatonic regions of voice-leading with the Boretz regions/spiders. Cohn also describes generalized parsimonious voice-leading models for chords with arbitrary cardinality $n$ in a tonal system with $m$ available tones in an octave such that $n$ is a divisor of $m$. Within the standard 12-tone system, he discusses in particular the existence of a voice-leading model similar to the Weitzmann/hexatonic and Boretz/octatonic systems but for the $n = 6$ case involving the nearly symmetric mystic and Wozzeck chords.

[1.3] In this paper, I develop a visual representation for voice-leading parsimony in the mystic-Wozzeck genus, constructing a unified model of dodecatonic cycles and centipede voice-leading regions. The mystic-Wozzeck genus and its voice-leading regions can be generated using a perturbative method similar to what Childs (1998) uses for seventh chords which Cohn (2012) repurposes for major and minor triads. In section II, I discuss symmetric partitions of the octave and perturbations of chords generated from such symmetries. In section III, I recapitulate Cohn's work on the generation of a unified model for voice-leading parsimony in the minor-major genus, starting with the symmetric chords from section II. I also consider available Neo-Riemannian transformations that are used to voice-lead within Weitzmann regions/waterbugs and within hexatonic cycles. Section IV extends the approach of the previous section to the $n = 4$ case, and I walk through the perturbative construction of the Boretz regions/spiders and the octatonic cycles as they are represented in Child's (1998) article and Cohn's (2012) \textit{Audacious Euphony}. I also propose a procedure for reducing the dimension of Childs' (1998) cubic diagram of the octatonic cycle.

[1.4] While the voice-leading regions for the $n = 3$ and $n = 4$ nearly symmetric chords are well-known, I walk through their derivations in this paper for the sake of highlighting the inherent similarities and self-consistencies between the $n = 3$, $n = 4$, and $n = 6$ cases. The hexatonic and octatonic voice-leading regions are relatively easy to visualize as in 2 and 3 dimensions, respectively. But, this is not the case with my proposed dodecatonic region, which would require 5 dimensions to represent the region as a convex polyhedron. My proposed method to reduce the dimension of Childs' octatonic region can be directly applied to the $n = 6$ case to allow for easy visualization of the dodecatonic region in 2 dimensions. Walking through the $n = 4$ case makes this dimension-reduction process especially transparent. The extension of the $n = 3$ and $n = 4$ derivations to the $n = 6$ case is carried out in Section V, forming the central portion of this paper. I construct what I refer to as \textit{centipede} regions and \textit{dodecatonic} cycles for voice-leading between mystic and Wozzeck chords. Section VI discusses the set-theoretic properties of the hexatonic, octatonic, and dodecatonic regions and proposes a future direction for research.
\begin{figure*}
  \centering
  \includegraphics[width=\textwidth]{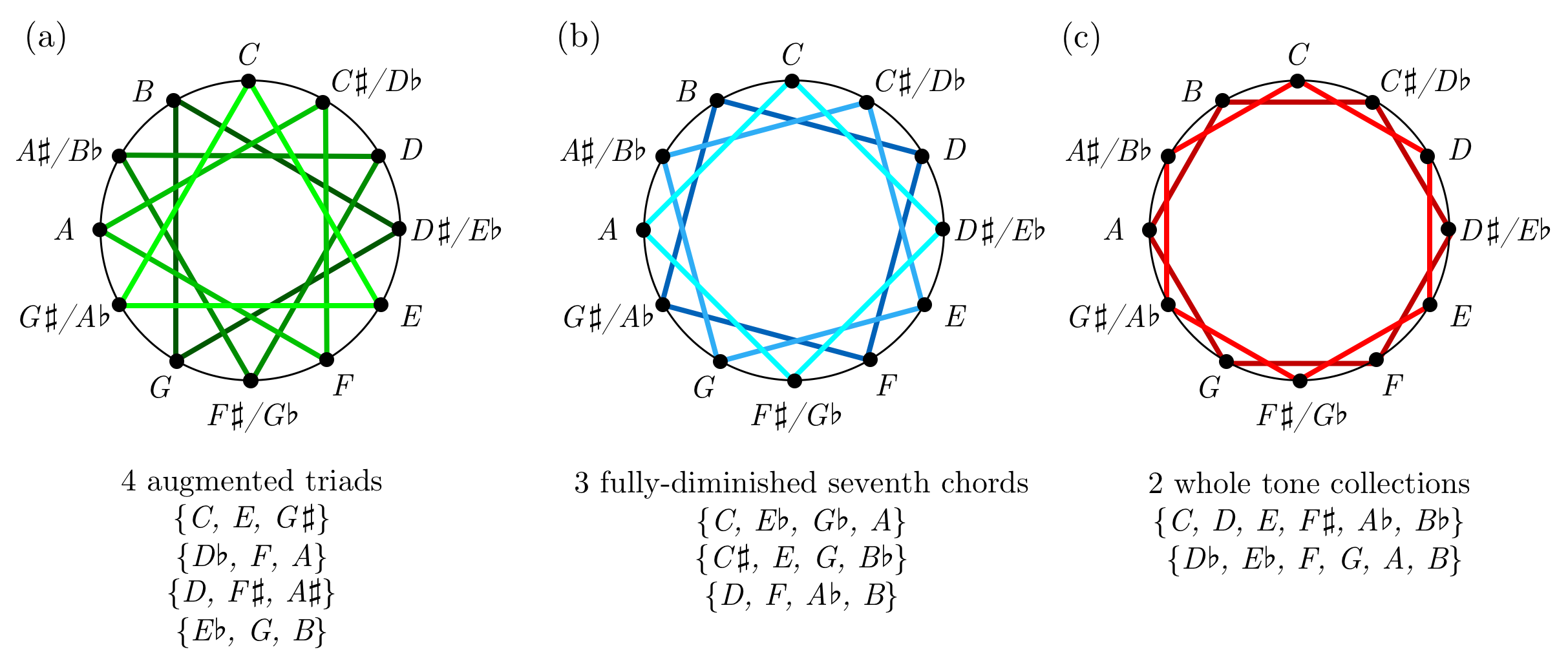}
  \caption{Possible symmetric partitions of the octave for the (a) $n =3$, (b) $n = 4$, and (c) $n = 6$ cases.}
\end{figure*}
\section*{II. Symmetric Partitions of the Octave}
[2.1] The splitting of an octave into 12 tones allows for the generation of interesting symmetries that often translate into musically relevant structures. For example, in mathematics, the cyclic group $\mathbb{Z}_{12}$ consists of the integers from 0 to 11, and elements of the group are related to each other by addition modulo 12. Of utmost relevance to music theory is the fact that $\mathbb{Z}_{12}$ has generators 1, 5, 7, and 11. Using the generator 1, one may construct all of the elements of the group by starting with one element---say, 0---and adding (modulo 12) the generator repeatedly:
\begin{align*}
0 &\equiv 0 \mod 12 \\
1 &\equiv (0 + 1)\mod 12 \\
2 &\equiv (1 + 1)\mod 12 \\
&\dots\\
11 &\equiv (10 + 1)\mod 12.
\end{align*}
If the generator is 7, the group can be constructed with the same technique, but the elements start to appear in a different order:
\begin{align*}
0 &\equiv 0 \mod 12 \\
7 &\equiv (0 + 7)\mod 12 \\
2 &\equiv (7 + 7)\mod 12 \\
9 &\equiv (2 + 7)\mod 12 \\
&\dots\\
5 &\equiv (10 + 7)\mod 12.
\end{align*}
As commonly done in musical set theory, if one assigns a pitch class to each integer (\textit{i.e.}, $C=0$, $C\sharp=1$, $\dots$, $B=11$), the order in which the elements are generated starts to form either the chromatic circle or the circle of fifths. Generators 5 and 11, likewise, form the circle of fourths and the descending chromatic circle, respectively. Hence, the mathematics can serve as a tool for the rigorous construction of well-known musical phenomena.

[2.2] The number 12 has divisors 1, 2, 3, 4, 6, and 12. The divisor $n$ indicates how the notes in the octave may be partitioned symmetrically. For the case where $n = 3$, one starts with a particular note---perhaps $C$---and selects every $12/3 = 4$th note that appears in the chromatic circle\footnote[2]{Note here that any circle generated by 1, 5, 7, and 11 may be used for this purpose.}. Connecting these notes with lines, one sees that an equilateral triangle is formed. Furthermore, there are 4 individual equilateral triangles that can be formed, and no vertices intersect. \textbf{Figure 1(a)} shows that symmetrically partitioning the octave for the $n = 3$ case yields the 4 augmented triads with nonintersecting sets of pitch classes: $\{C, E, G\sharp\}$, $\{D\flat, F, A\}$, $\{D, F\sharp, A\sharp\}$, and $\{E\flat, G, B\}$.
\begin{figure*}
  \centering
  \includegraphics[width=\textwidth]{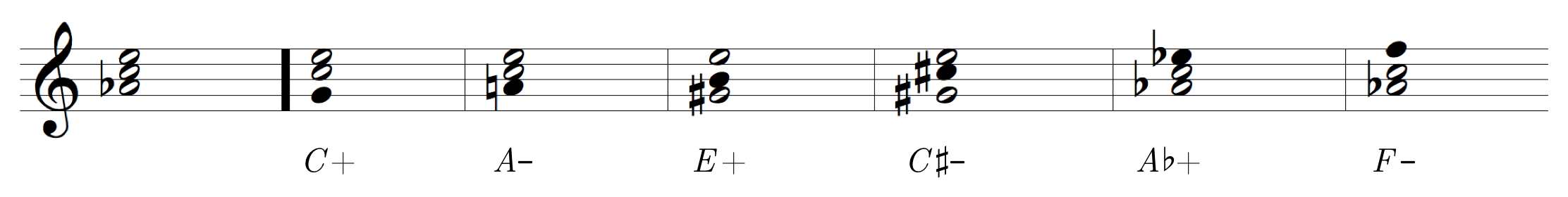}
  \caption{Perturbations of the $A\flat$ augmented triad, generating a Weitzmann region of major and minor chords.}
\end{figure*}

[2.3] For the $n = 4$ case, one finds the symmetric chords in the same manner. This time, the cardinality of the chord is also $n = 4$, and one must select every $12/4 = 3$rd note in the chromatic circle to complete the chord, which is a fully diminished seventh chord. \textbf{Figure 1(b)} clearly illustrates that connecting the notes in a chord generates a square, and there are 3 independent squares that share no pitch classes with each other. These chords are comprised of the collections $\{C, E\flat, G\flat, A\}$, $\{C\sharp, E, G, B\flat\}$, and $\{D, F, A\flat, B\}$.

[2.4] The $n = 6$ case is determined using the same method. By choosing every other note ($12/6 = 2$), one constructs a whole-tone scale, which geometrically forms a regular hexagon superimposed onto the chromatic circle. \textbf{Figure 1(c)} shows the 2 independent whole-tone scales, which are comprised of $\{C, D, E, F\sharp, A\flat, B\flat\}$ and $\{D\flat, E\flat, F, G, A, B\}$.
\begin{figure*}
  \centering
  \includegraphics[width=0.5\textwidth]{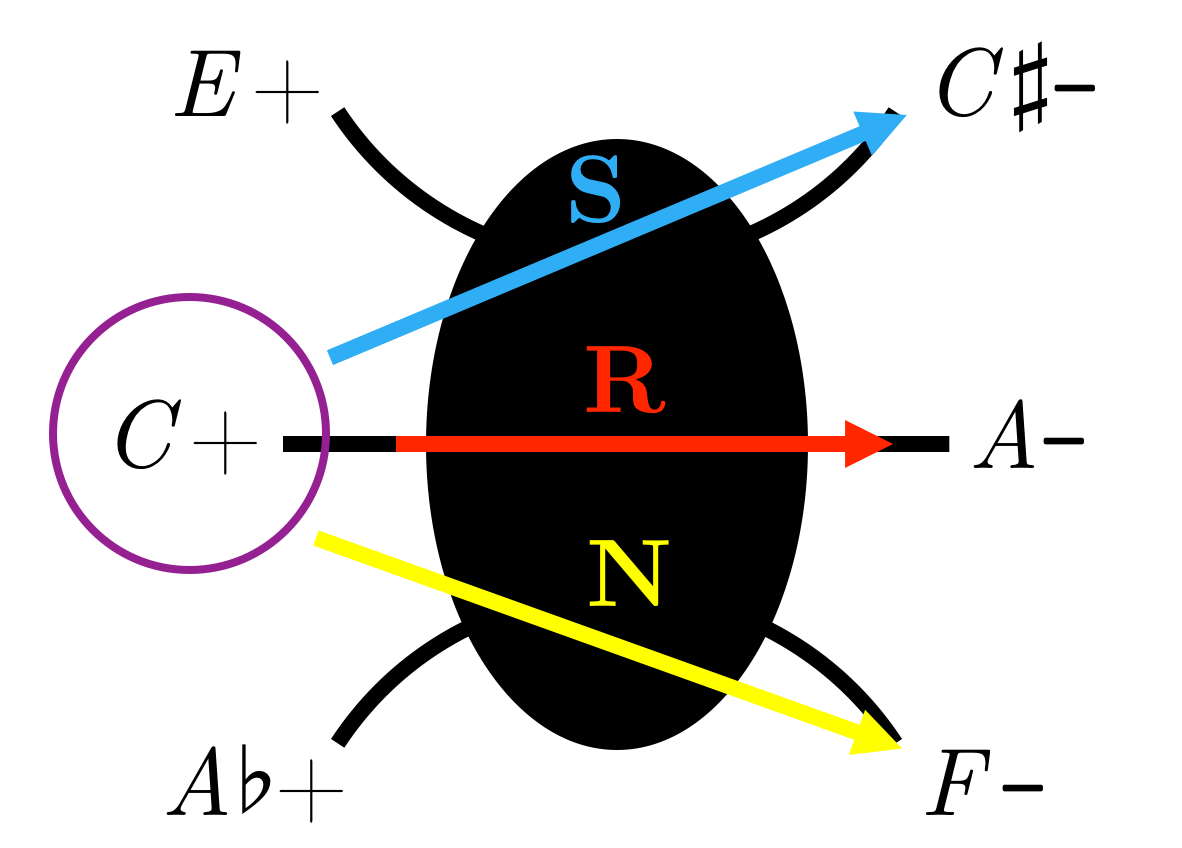}
  \caption{A Weitzmann waterbug, adapted from Cohn (2012).}
\end{figure*}
\begin{figure*}
  \centering
  \includegraphics[width=0.9\textwidth]{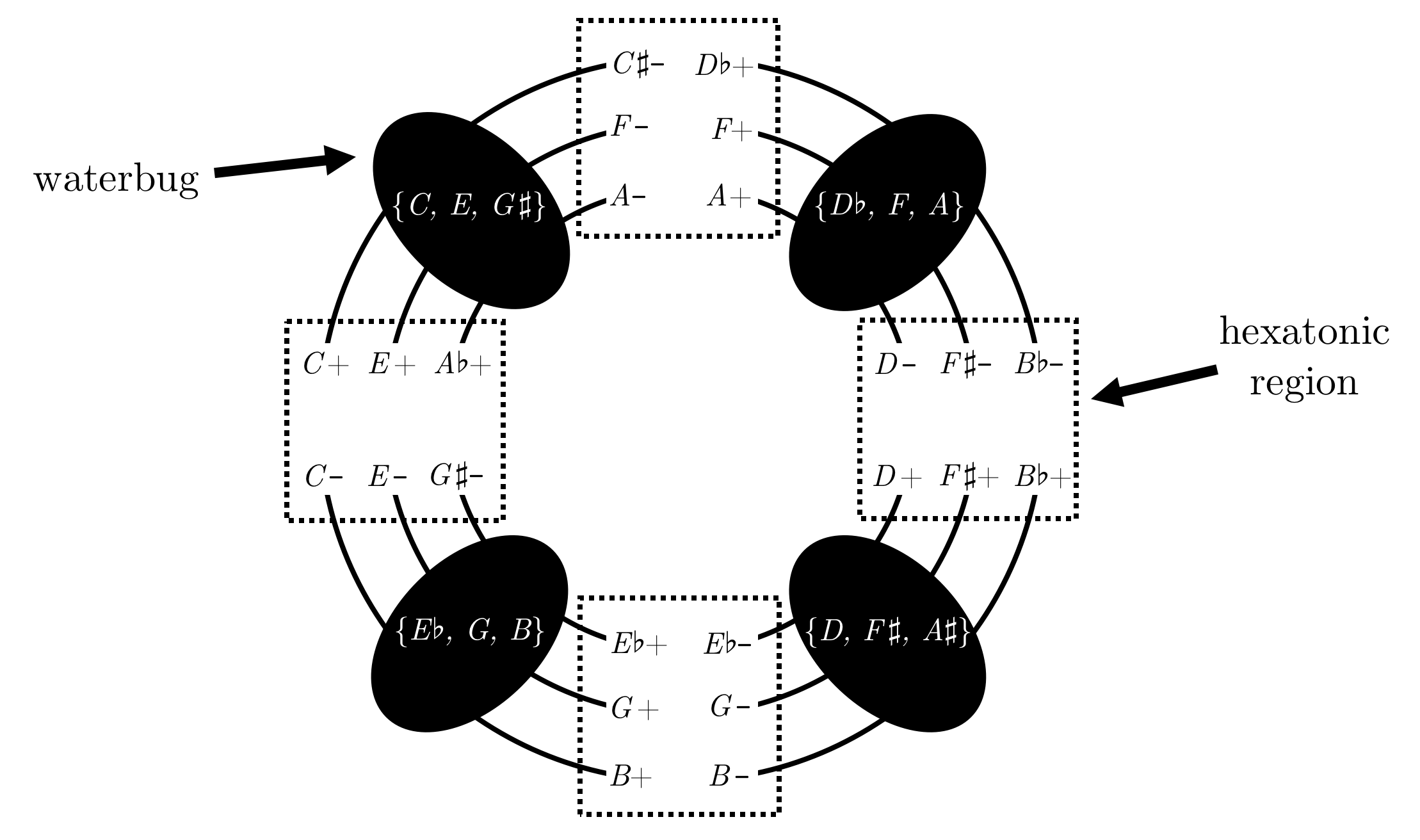}
  \caption{Reproduction of Cohn's (2012) unified voice-leading model for nearly symmetric triads.}
\end{figure*}

[2.5] The $n = 1$, 2, and 12 cases function similarly, but voice-leading between classes of single note, dyads, and nearly-full chromatic sets is relatively uninteresting.

\subsection*{Nearly Symmetric Chords}
[2.6] The correspondence between these symmetric geometries and well-known pitch collections demonstrates the musical relevance of certain mathematical structures. Each of the following three sections begins with a preliminary discussion on perturbation of the symmetric chords previously discussed. As Dmitri Tymoczko (2011) heavily emphasizes in his book \textit{A Geometry of Music}, many of fundamental triads, seventh chords, and scales commonly used Western music are nearly even chords; that is to say, they are not quite symmetric chords, but they are only a few semitone displacements away from symmetric chords. While Tymoczko has a relatively open definition of ``nearly even,'' in this paper I specifically focus on chords that are exactly a single-semitone displacement (SSD) from one of the symmetric chords previously discussed. I shall refer to this specific class of perturbed chords as \textit{nearly symmetric}. Nearly symmetric chords of a given cardinality $n$ demonstrate a specific, consistent pattern of voice-leading.

\section*{III. Weitzmann Waterbugs and Hexatonic Cycles}
[3.1] The literature on parsimonious voice-leading between major and minor triads is extensive, as the field of Neo-Riemannian theory essentially developed from this central topic. Richard Cohn's (1996) article on hexatonic cycles thoroughly develops this topic, focusing on the set of maximally smooth cycles that are generated from single-semitone displacements of the major and minor triads. For the collection of 24 major and minor triads, 4 independent maximally smooth cycles---called hexatonic cycles---each consisting of 3 major and 3 minor chords can be constructed via a simple procedure: one starts with a consonant triad and needs to voice lead to another consonant triad only utilizing single-semitone displacements. This results in a transformation that is more specifically an involution: a major triad can be transformed to only a minor triad via a single-semitone displacement, and a minor can be transformed to only a major triad via a single-semitone displacement. This procedure can be repeated until one returns to the starting chord. A maximally smooth cycle of major and minor triads consists of exactly 6 chords, and there are 4 independent maximally smooth cycles for this collection.

[3.2] There are multiple ways in which the elements of the hexatonic cycle can be constructed; the procedure presented above is one of them. I now describe a procedure that Cohn (2012) models after Childs' (1998) approach: namely, I exploit the near symmetric nature of the major and minor chords to derive the voice-leading regions.

[3.3] In the previous section, I showed that for the $n = 3$ case, there are 4 independent augmented triads that symmetrically partition the octave. \textbf{Figure 2} shows the triad $\{A\flat,C,E\}$. If the $A\flat$ is perturbed down a semitone, a $C$ major triad is obtained, and it is written as $C+$. If the $A\flat$ is perturbed up a semitone, an $A$ minor triad is obtained, notated as $A-$. Performing downward and upward perturbations on the pitch $C$ results in $E+$ and $C\sharp-$, and perturbing the pitch $E$ generates $A\flat+$ and $F-$ from the augmented triad. This collection of chords $\{C+,A-,E+,C\sharp-,A\flat+,F-\}$ is known as the Weitzmann region (Cohn 2012). The 3 other independent augmented triads similarly generate 3 separate Weitzmann regions, and no two Weitzmann regions share any chords.

[3.4] An well-known visual representation of the Weitzmann region is the Weitzmann waterbug, shown in \textbf{Figure 3} (Cohn 2012). The legs on one side of the waterbug's body correspond to the $(+)$ chords, while the legs on the other half correspond to the $(-)$ chords. The familiar Neo-Riemannian transformations that act within this region are $\vb{R}$ (relative), $\vb{N}$ (\textit{Nebenverwandt}), and $\vb{S}$ (slide).  Applying $\vb{R}$ to $C+$ requires movement of a single voice by 2 semitones to produce $A-$. Using Douthett and Steinbach's (1998) formal definition of $P_{m,n}$-related chords, one would say that $C+$ and $A-$ are $P_{0,1}$ related because one can shift between the two chords only by moving 1 note by a whole tone and 0 voices by semitone. The transformations $\vb{N}$ and $\vb{S}$ both require moving two voices in parallel motion by a single semitone (each). $\vb{S}$ moves the perfect 5th, shifting the pitches $C$ and $G$ to $C\sharp$ and $G\sharp$ (so $C+$ becomes $C\sharp-$), while the pitch $E$ remains invariant. $\vb{N}$ shifts the minor 3rd, $E$ and $G$ are transformed to $F$ and $A\flat$ (so $C+$ becomes $F-$), leaving the root invariant. Thus, chords related by $\vb{N}$ or $\vb{S}$ are said to be $P_{2,0}$ related.

[3.5] The four Weitzmann waterbugs that are generated from the 4 augmented triads have no intersection, so to represent the relationships between the full system of waterbugs, one can use a diagram like Douthett and Steinbach's Cube Dance (1998) or Cohn's (2012) unified waterbug/hexatonic figure in \textit{Audacious Euphony}. For the sake of visual clarity, only Cohn's diagram is recreated here in \textbf{Figure 4}, even though Cube Dance contains additional voice-leading information. Cohn's diagram is constructed simply by placing the waterbugs in a square such that the root names match at the ``bridge'' regions of the waterbugs. For example, the left bridge region contains the chords $\{C+, E+, G+\}$ from one waterbug and $\{C-, E-, G-\}$ from the adjacent waterbug. Arranging the waterbugs in this way reveals that the chords within the bridge regions are exactly the ones in the hexatonic cycles (Cohn 2012). Moreover, all 24 major and minor chords are represented in this figure, so \textbf{Figure 4} is indeed a unified model of voice-leading for the major/minor collection.
\begin{figure*}[t]
  \centering
  \includegraphics[width=\textwidth]{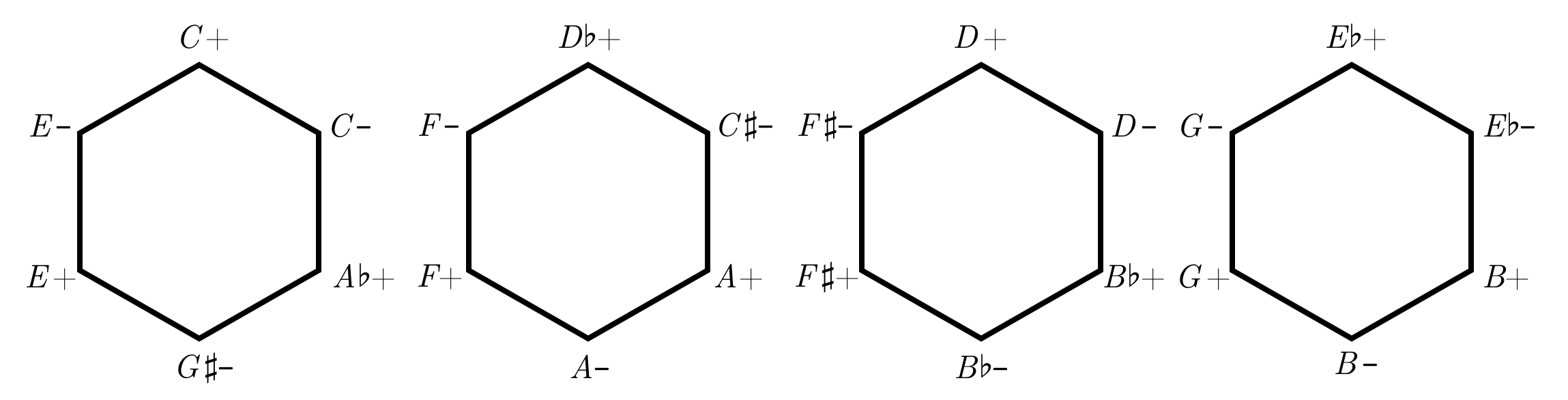}
  \caption{Hexatonic cycles shown as hexagons (Cohn 2000).}
\end{figure*}
\begin{figure*}[t]
  \centering
  \includegraphics[width=0.3\textwidth]{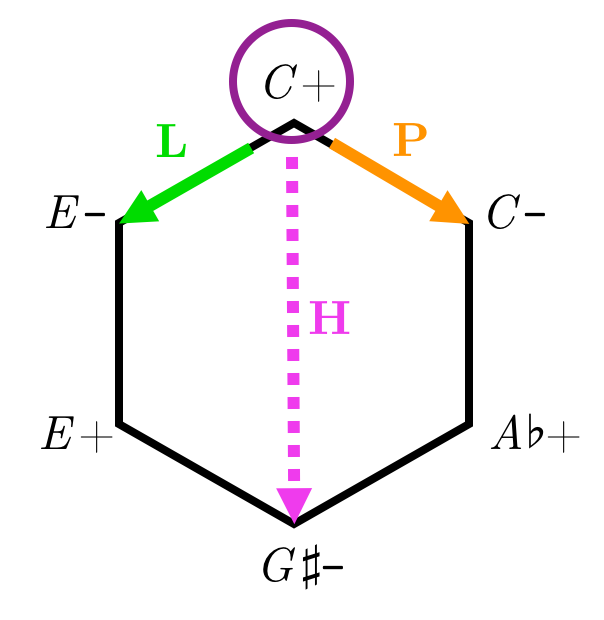}
  \caption{Available voice-leading transformations in a hexatonic region (Cohn 1996).}
\end{figure*}

[3.6] The hexatonic regions (or cycles) themselves are constructed from the familiar Neo-Riemannian transformations. The four hexatonic cycles, which Cohn (1996; 2000) names ``Northern,'' ``Southern,'' ``Eastern,'' and ``Western,'' are shown as hexagons instead of circles in \textbf{Figure 5}. In his 1996 paper, Cohn draws these hexatonic regions as circles and the 2000 article presents them as hexagons. For the sake of geometric consistency with the $n = 4$ and $n = 6$ cases which I will present, I choose the hexagonal representation. This is particularly useful because a ``musical'' meaning can be attributed to each vertex and edge of the hexagon: every vertex corresponds to a major or minor chord, and each straight line (forming the sides of the hexagon) corresponds to the identical voice leading distance. One can say that two chords $X$ and $Y$ connected by a straight line in a given hexatonic region are $P_{1,0}$-related, since one can construct $Y$ from $X$ (and $X$ from $Y$) by only moving 1 note by a semitone. The transformations $\vb{P}$ (parallel) and $\vb{L}$ (\textit{Leittonweschel}) are responsible for transformations between $P_{1,0}$-related chords. Starting with a chord of a given modality, the only chord in the hexatonic region of the opposite modality that is not $P_{1,0}$-related to the starting chord is the hexatonic pole, reached by the $\vb{H}$ transformation. The geometric functions of these 3 transformations on a sample chord, $C+$, are shown in \textbf{Figure 6}. I will show that in the next two sections, construction of the voice-leading regions for the $n = 4$ and $n = 6$ cases follow the same procedure as that for the $n = 3$ case, and the musical interpretation of the geometric elements (e.g. straight lines and $P_{m,n}$-relatedness) remains consistent.
\begin{figure*}[t]
  \centering
  \includegraphics[width=\textwidth]{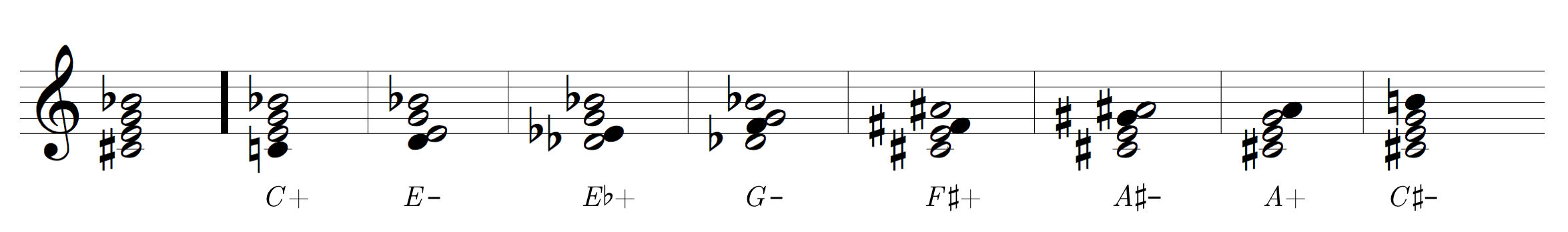}
  \caption{Perturbations of the $F\sharp$ fully diminished seventh chord, generating a Boretz region of dominant seventh and half-diminished seventh chords.}
\end{figure*}
\begin{figure*}[t]
  \centering
  \includegraphics[width=0.5\textwidth]{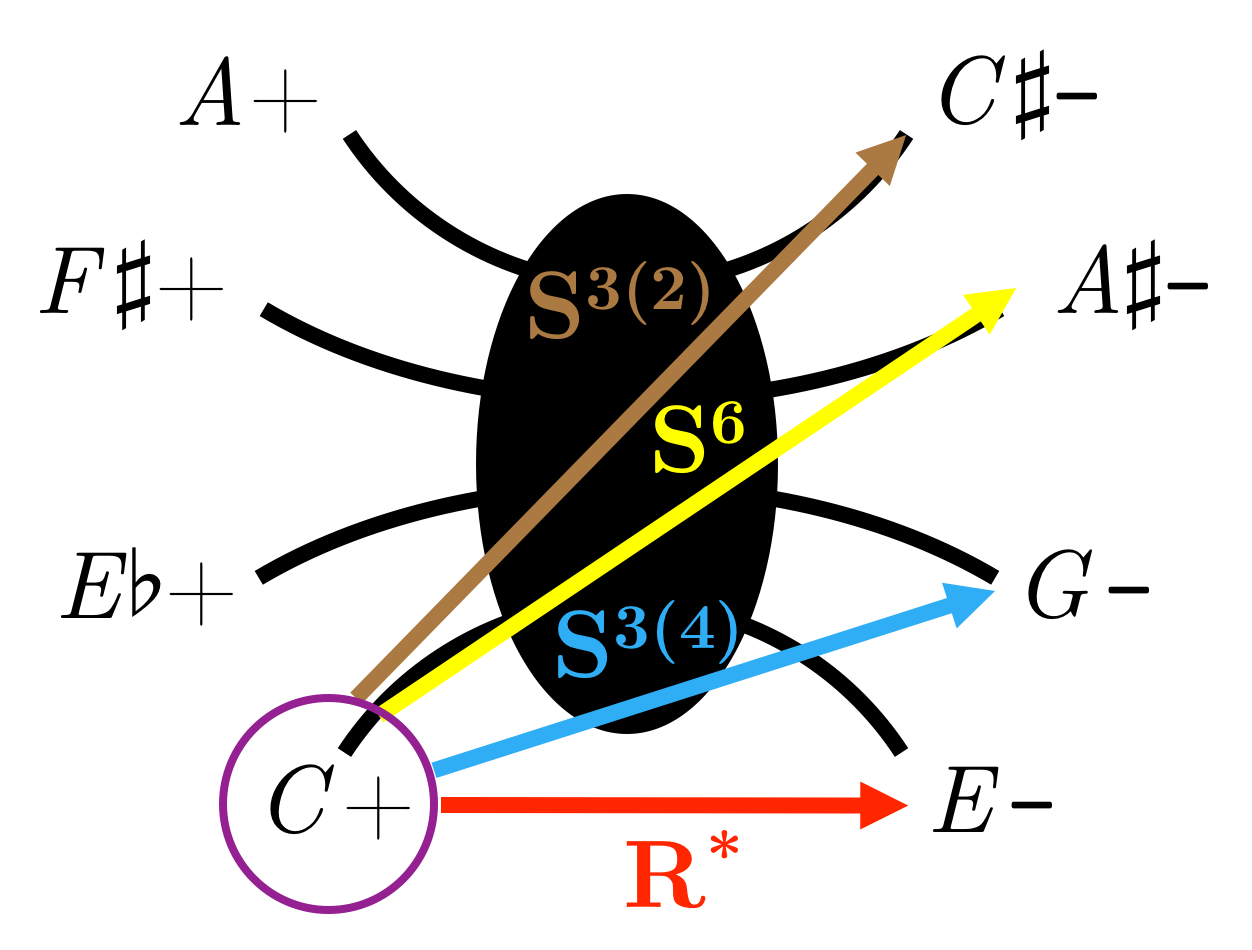}
  \caption{A Boretz spider (Cohn 2012).}
\end{figure*}
\section*{IV. Boretz Spiders and Octatonic Cycles}
[4.1] Adrian Childs' (1998) paper extends the previously described unified voice-leading model to the $n = 4$ case. In this section, I arrive at Childs' results using the same method outlined in the previous section, starting from symmetric partitions of the octave\footnote[3]{In fact, Childs first proposed the perturbative approach used to generate the nearly-symmetric seventh chords in his 1998 paper. Cohn adapted this approach to the triadic case, which he presented in the first edition of \textit{Audacious Euphony}, which was published in 2000. In this paper, though I present the derivation of the $n = 3$ voice-leading regions first, the perturbative approach to the $n = 4$ indeed was published first.}. In section II, I showed that for the $n = 4$ case, the octave is symmetrically partitioned independently by 3 fully diminished chords. The perturbations of one of these possible fully diminished chords are shown in \textbf{Figure 7}. Notationally, the symbols $(+)$ and $(-)$ are used to represent the dominant seventh and half-diminished seventh chords, respectively.
\begin{figure*}[t]
  \centering
  \includegraphics[width=0.8\textwidth]{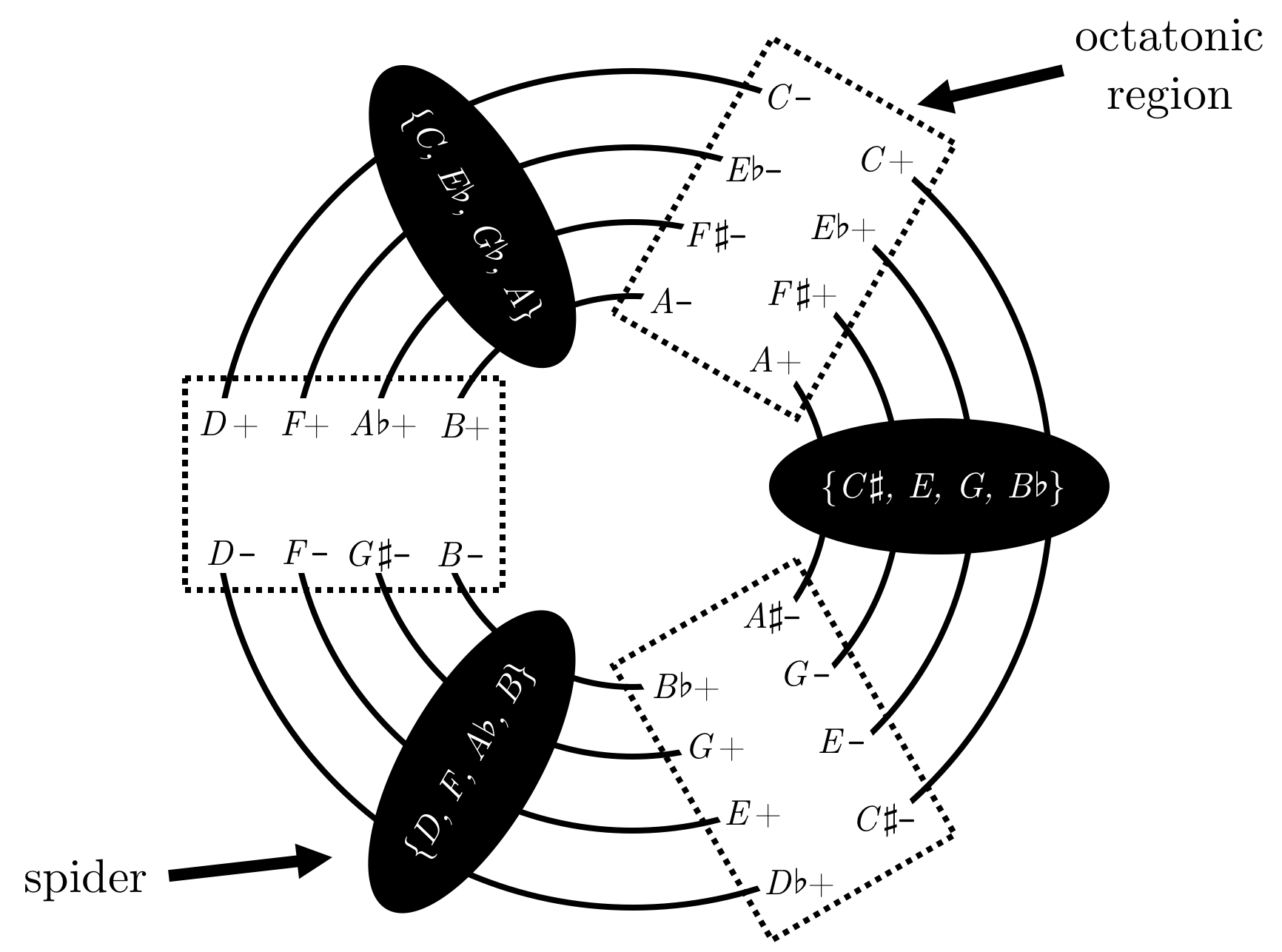}
  \caption{Recreation of Cohn's (2012) unified voice-leading model for nearly symmetric seventh chords.}
\end{figure*}

[4.2] Starting with $C\sharp$ fully diminished seventh, if the bottom $C\sharp$ is perturbed downward, the resulting chord is $C$ dominant seventh, written as $C+$. Perturbing the $C\sharp$ upward results in $E$ half-diminished seventh, written as $E-$. Perturbing the $E$ downward results in $E\flat+$, and shifting the same note upward results in $G-$. In general, a downward perturbation results in a dominant seventh chord, and an upward perturbation results in a half-diminished seventh chord. Continuing these perturbations for the remaining two notes in the initial $C\sharp$ fully diminished seventh chord, one can generate a Boretz region, given by the collection of chords $\{C+, E-, E\flat+, G-, F\sharp+, A\sharp-, A+, C-\}$ (Cohn 2012). There are 3 independent Boretz regions, and the union of these 3 sets gives the full collection of all dominant seventh and half-diminished seventh chords.

[4.3] As the Weitzmann waterbug conveniently represents the Weitzmann region, the Boretz spider, shown in \textbf{Figure 8}, is the $n = 4$ visual representation of the Boretz region (Cohn 2012). In each of the 3 Boretz spiders, a dominant seventh chord corresponds to a leg on one half of the spider's body, and the half-diminished seventh chords are placed on the other half of the body. Well-defined involutions can be applied to any given chord in order to transform one chord to any other chord on the other side of the Boretz spider's body. These particular transformations are $\vb{R}^*$, $\vb{S^{3(4)}}$, $\vb{S^{6}}$, and $\vb{S^{3(2)}}$ (Childs 1998; Cohn 2012).

[4.4] The $\vb{R}^*$ transformation is the $n = 4$ analogy to the triadic $\vb{R}$ transformation. Though there is no formal definition of a ``relative'' half-diminished seventh chord for a given dominant seventh chord or vice versa, the $\vb{R^*}$ transformation moves 1 voice by 2 semitones, just like $\vb{R}$. Thus, $C+$ and $E-$ are $P_{0,1}$-related.

[4.5] The remaining $\vb{S}$-type transformations for the Boretz spiders have functions that are analogous to the triadic $\vb{S}$ and $\vb{N}$ transformations. In $\vb{S^{3(4)}}$, $\vb{S^{6}}$, and $\vb{S^{3(2)}}$, one ``slides'' 2 voices by 1 semitone in parallel motion. Thus, chords that are related by $\vb{S^{3(4)}}$, $\vb{S^{6}}$, and $\vb{S^{3(2)}}$ are said to be $P_{2,0}$-related, just like in the $n = 3$ case. The first number in the superscript (\textit{e.g.} the ``3'' in $\vb{S^{3(4)}}$) refers to the interval within the 4-note chord that is held invariant. The number in parentheses denotes the interval that ``slides'' due to the transformation. For example, if one applies $\vb{S^{3(4)}}$ to $F+$, then the set-theoretic interval 3, which is a minor 3rd, is held invariant. There are 2 minor 3rds---pitch class $E$ to $G$ and $G$ to $B\flat$---so one looks at the number $4$ to determine the interval that is shifted. The interval $4$ is a major 3rd, so pitch classes $C$ and $E$ must be shifted. A downward shift does not result in a dominant seventh or half-diminished seventh chord, so $\vb{S^{3(4)}}$ specifically transforms $C+$ to $G-$ and vice versa. As a comparison, $\vb{S^{3(2)}}$ leaves the minor 3rd $E$ to $G$ invariant while shifting $D$ and $B\flat$, transforming $C+$ to $C\sharp-$. Cohn (2012) writes $\vb{S^{6(5)}}$ as $\vb{S^6}$ because it is implied that invariance of the interval 6 requires shifting of the perfect 4th, as there is only 1 possible tritone within a dominant seventh or half-diminsihed seventh chord.
\begin{figure*}[t]
  \centering
  \includegraphics[width=0.8\textwidth]{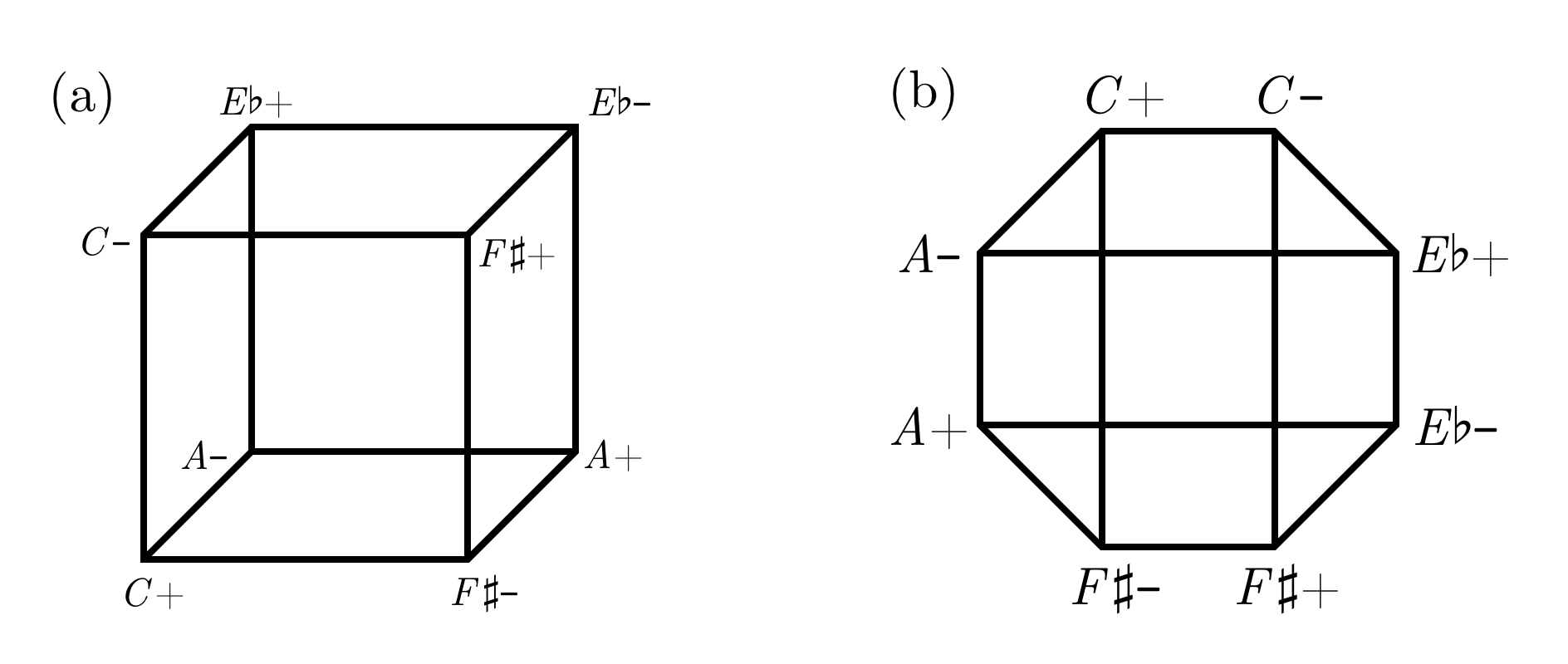}
  \caption{A octatonic region represented as a (a) cube (Childs 1998) and (b) a 2-dimensional network.}
\end{figure*}
\begin{figure*}[t]
  \centering
  \includegraphics[width=0.9\textwidth]{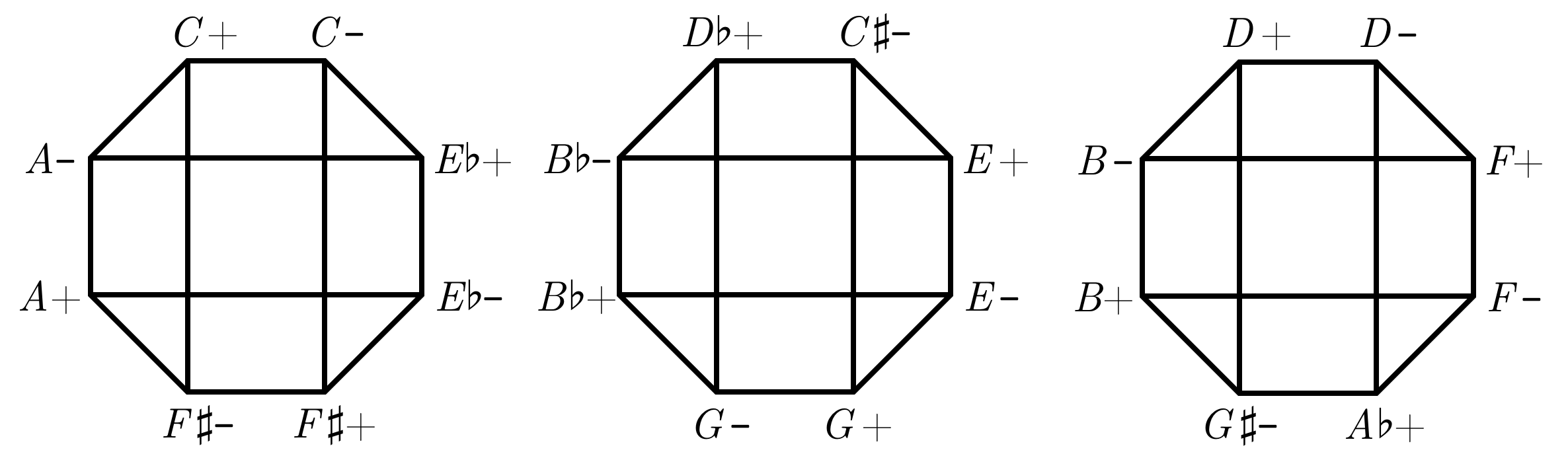}
  \caption{Octatonic cycles (Childs 1998), shown as 2-dimensional graphs.}
\end{figure*}

[4.6] Cohn (2012) joins the 3 independent Boretz spiders in a unified voice-leading model for the $n = 4$ case, and Douthett and Steinbach (1998) have a similar figure---Power Towers---which shows additional voice-leading capabilities within the bridge region. For visual simplicity, Cohn's figure is reproduced in \textbf{Figure 9}. As with the triads, the $(+)$ chord root names on one spider must be ``bridged'' with the $(-)$ chords with the same root names on another spider. For example, the bottom bridge region of \textbf{Figure 9} unites $\{C\sharp+, E+, G+, B\flat+\}$ on one spider with $\{C\sharp-, E-, G-, A\sharp-\}$ on the adjacent spider. The 3 bridge regions that are generated are---in analogy with the $n = 3$ case---referred to as octatonic regions. Childs (1998) shows that each octatonic region can be displayed as a cubic network, and one of the three cubes is constructed in \textbf{Figure 10(a)}. Geometrically, a cube is a desirable structure for describing this type of chord collection because it is not only a convex polytope, but it is also possible to assign its vertices to $(+)$ and $(-)$ chords such that no $(+)$ vertex has any edge shared with another $(+)$ vertex, and no $(-)$ vertex has any edge shared with another $(-)$ vertex.

[4.7] Although a cube is convenient for visualizing the $n = 4$ voice-leading region, visualizing convex polyhedra in $n > 3$ spatial dimensions becomes an impossible task. A method for reducing the dimension of higher-dimensional geometric structures will surely prove useful when trying to visualize the $n = 6$ case, since 5 spatial dimensions would be required. Thus, I propose a method for flattening such geometric representations of bridge regions to 2-dimensions; this becomes an especially powerful and useful tool when dealing with the $n = 6$ case, and it can also be applied to $n = 4$. (The $n = 3$ hexatonic cycles already have a 2-dimensional representation.) The geometric structure shown in \textbf{Figure 10(b)} is my alternative to Childs' (1998) cubic network for the octatonic region. Since there are 8 chords in the octatonic collection, an octagon provides the neatest ``frame'' for the structure, just like a hexagon does for the $n = 3$ case. I arrange the $(+)$ and $(-)$ chords around the octagon in chromatically sequential order, alternating between the $(+)$ and $(-)$ chords, just like Cohn does for the hexatonic cycles. However, unlike the hexatonic cycles, each vertex in the octagon is connected to more than 2 other vertices. Examining Childs' cubic network, I draw lines \textit{inside} the perimeter of the octagon to connect chords that also have connections in the cube. As a result, any $(+)$ chord is connected to all other $(-)$ chords except its own octatonic pole, and any given $(-)$ chord is connected to all other $(+)$ chords except its own octatonic pole. The 3 independent octatonic regions are shown in the 2-dimensional form in \textbf{Figure 11}.
\begin{figure*}[t]
  \centering
  \includegraphics[width=0.4\textwidth]{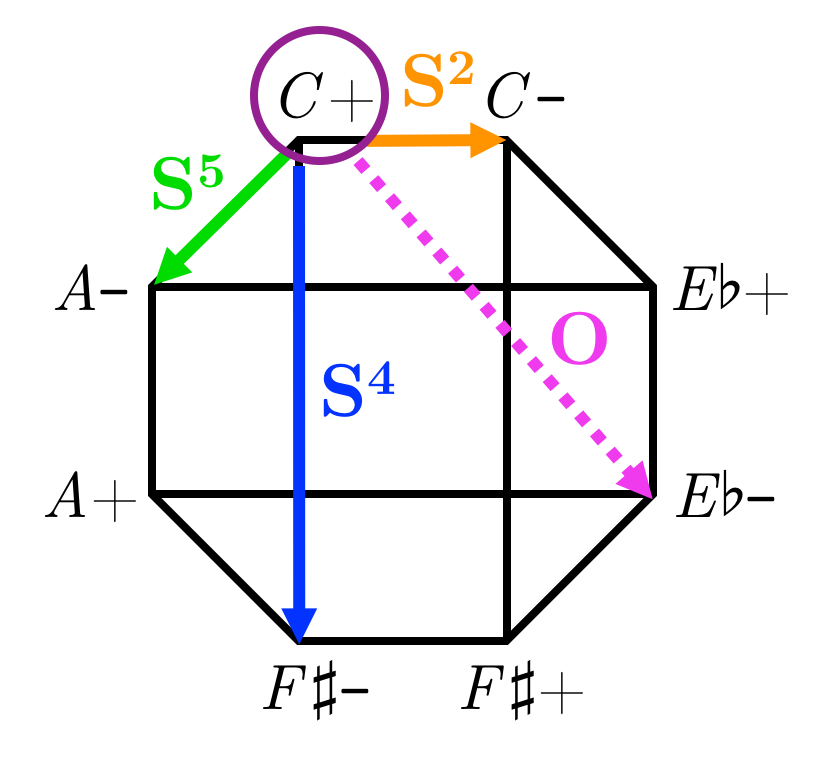}
  \caption{Available voice-leading transformations in an octatonic region (Childs 1998), shown in the 2-dimensional representation.}
\end{figure*}
\begin{figure*}[t]
  \centering
  \includegraphics[width=\textwidth]{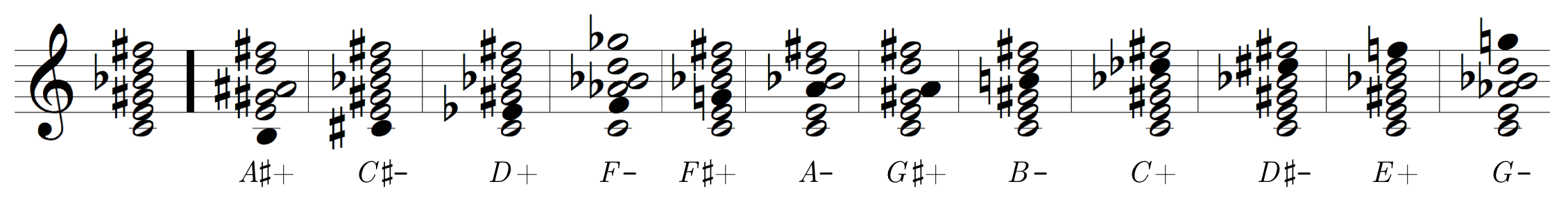}
  \caption{Perturbations of the whole-tone scale including $C$, generating a collection of Wozzeck and mystic chords.}
\end{figure*}

[4.8] The Neo-Riemannian transformations that act within this octatonic region are shown in \textbf{Figure 12}. Of the 4 involutions that are allowed, 3 are ``slide'' transformations: $\vb{S^2}$, $\vb{S^4}$, and $\vb{S^5}$. These 3 chords follow the same notational scheme as the $\vb{S}$-type transformations that acted within the Boretz region, and they are abbreviations of $\vb{S^{2(3)}}$, $\vb{S^{4(3)}}$, and $\vb{S^{5(6)}}$. The $\vb{S^2}$ transformation holds the major 2nd (set-theoretic interval 2) invariant while sliding the minor 3rd, so the chord $C+$ would be transformed to $C-$. $\vb{S^4}$ would hold the major 3rd invariant while shifting the minor 3rd, so $C+$ transforms to $F\sharp-$. Likewise, $\vb{S^5}$ would transform $C+$ to $A-$, since the perfect 4th is held invariant while the tritone is shifted. In each of the above cases, the chords related by these $\vb{S}$-type transformations are $P_{2,0}$-related since two voices are shifted by 1 semitone each, and 0 voices are shifted by a whole tone.

[4.9] Within the 2 geometric representations of the octatonic region in \textbf{Figure 10}, a solid line between two chords represents an identical voice-leading distance. That is to say, any two chords connected by a solid line are $P_{2,0}$-related. This fact may seem obvious from Childs' (1998) cubic diagram; but in the 2-dimensional reduction, the correspondence between length of a connecting line and voice-leading distance is lost. While this sacrifice of information must be made in the reduction of dimension, the flattening of geometric structures to 2 dimensions provides a powerful and, most importantly, geometrically consistent method for treating the $n = 3$, $n = 4$, and---as I will show in the next section--- $n = 6$ case on equal footing, visually. 

[4.10] Nonetheless, it is easy to see that given a particular chord in an octatonic region, there is only one chord of the opposite modality that is not $P_{2,0}$-related to the starting chord. This is known as the octatonic pole (Childs 1998), and I refer to it as the $\vb{O}$ transformation\footnote[7]{Neither Childs (1998) nor Cohn (2012) explicitly associated a letter with the octatonic pole transformation. In analogy with $\vb{H}$ for hexatonic pole, I coin $\vb{O}$ for octatonic pole.}. On Childs' cubic representation, the octatonic pole is present on the vertex that is farthest away from the starting chord. In the 2-dimensional representation, the $\vb{O}$ transformation connects two chords of the opposite modality that are not joined by a solid line. Like the hexatonic pole for triads, a chord and its octatonic pole share no pitch classes. $\vb{O}$ applied to $C+$ results in $E\flat-$.

[4.11] Discussion of an octatonic region also motivates a search for maximally smooth cycles. In a single hexatonic region, there is only one allowed hexatonic cycle: the alternating $\vb{P}$-$\vb{L}$ cycle. Hexatonic cycles are aptly named because all of the chords in a hexatonic region are utilized in a maximally smooth cycle, a definition which was formalized by Cohn (1996). In the octatonic region, one will notice that any closed path with no loops is indeed a maximally smooth cycle, according to Cohn's definition. The union of all pitch classes contained in the chords of any one of these maximally smooth cycles is an octatonic collection given by set class 8-28, so these maximally smooth cycles can be called octatonic cycles.

\section*{V. Centipedes and Dodecatonic Cycles}
[5.1] This section forms the crux of this paper. Here, I present the perturbative derivation of nearly symmetric hexachords that comprise the mystic-Wozzeck genus. These chords exhibit voice-leading properties similar to the major/minor triads and dominant/half-diminished seventh chords, and their voice-leading ``arthropod'' and bridge regions can be visually represented just as in the $n = 3$ and $n = 4$ cases.
\begin{figure*}[t]
  \centering
  \includegraphics[width=0.45\textwidth]{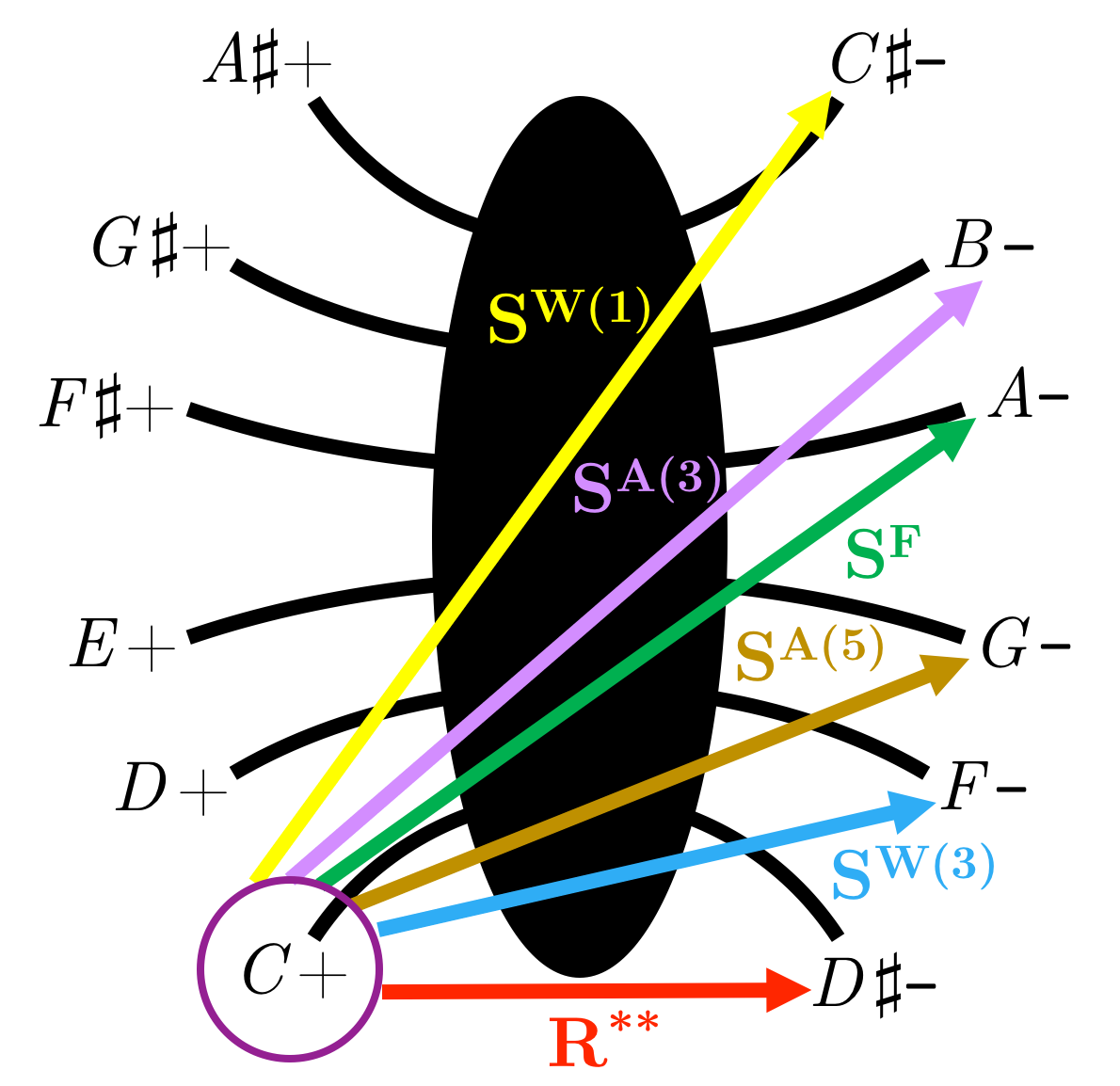}
  \caption{A centipede generated from perturbations of a whole-tone scale.}
\end{figure*}

\indent [5.2] Firstly, I must return to the symmetric perturbation of the octave discussed in section II. In the $n = 6$ case, there are only two ways to partition the octave, and these are the two non-intersecting whole-tone scales. \textbf{Figure 13} shows the perturbations of one of these whole-tone scales, given by the pitch class collection $\{C, D, E, F\sharp, G\sharp, B\flat\}$. As with the triads and seventh chords, if one perturbs any note in the whole-tone scale downward by a half step, the resulting chord is assigned a modality, and it is given the $(+)$ symbol. In this case, a downward perturbation of a note in the whole-tone scale results in a Wozzeck chord, where the name is taken from Alan Berg's opera \textit{Wozzeck} (Cohn 2012). An upward perturbation of any note in the whole-tone scale results in a mystic chord, which was brought to the attention of theorists primarily due to Alexander Scriabin's use of the chord in his compositions. The upward perturbation results in the opposite modality, so I notate mystic chords with the $(-)$ symbol.

[5.3] I assign an arbitrary naming scheme for the ``root'' of a mystic or Wozzeck chord so that the chords can be discussed simply by naming a letter name and modality symbol. In general, the root is the lower of the two notes in the minor 2nd interval in the mystic or Wozzeck chord. For example, the Wozzeck chord $\{C, D\flat, E, F\sharp, G\sharp, B\flat\}$ is notated as $C+$, and the mystic chord $\{C, D\flat, E\flat, F, G, A\}$ is notated as $C-$.

[5.4] Returning to \textbf{Figure 13}, one sees that individually perturbing the six notes of whole-tone scale both upward and downward results in six mystic chords and six Wozzeck chords. For the whole-tone scale beginning on $C$, the collection of nearly symmetric hexachords chords is given by $\{A\sharp+, C\sharp-, D+, F-, F\sharp+, A-,$ $G\sharp+, B-, C+, D\sharp-, E+, G-\}$. I propose that this collection can be represented visually as a centipede\footnote[9]{True centipedes do not typically have 12 legs, but a newborn garden symphylan (\textit{Scutigerella immaculata}) reportedly is indeed born with 6 pairs of legs but grows more over the course of its lifetime (Michelbacher 1938). The garden symphylan is commonly referred to as the garden centipede, which is why I have chosen ``centipede'' as a name for this voice-leading region.}, in analogy with the Weitzmann waterbug and Boretz spider.

[5.5] The ``centipede'' for mystic and Wozzeck chords is shown in \textbf{Figure 14}. As with the triads and seventh chords, the legs on one half of the centipede's body are all of the $(+)$ modality, and the $(-)$ chords are assigned to the legs on the other side of the body. It is easy to see that, given a starting chord, 5 of the 6 chords of the opposite modality are $P_{2,0}$-related to the starting chord, and the remaining chord of the opposite modality is $P_{0,1}$-related to the starting chord. This is directly analogous to the $n = 3$ and $n = 4$ cases. Thus, I can define Neo-Riemannian transformations that act on chords in the centipede that are analogous to the ``relative'' and ``slide'' transformations that act within the Weitzmann waterbug and Boretz spider.

[5.6] I define the $\vb{R^{**}}$ transformation as the ``relative'' transformation that connects two $P_{0,1}$-related chords of opposite modality. As with dominant seventh/half-diminished seventh chords, there is no formal definition of a ``relative'' mystic and Wozzeck chord, but the action of $\vb{R^{**}}$ is nonetheless well-defined: 1 voice must be moved by 1 step. As an example, $C+$ is transformed to $D\sharp-$ by $\vb{R^{**}}$.

[5.7] The 5 remaining Neo-Riemannian transformations in this region are ``slide'' transformations that involute between $P_{2,0}$-related chords: 2 voices are shifted in parallel motion by 1 semitone each. This means that 4 voices are invariant in the transformation. I define the following transformations following the notational convention of slide transformations in the $n = 4$ case: $\vb{S^{W(1)}}$, $\vb{S^{A(3)}}$, $\vb{S^{F}}$, $\vb{S^{A(5)}}$, and $\vb{S^{W(3)}}$. The following abbreviations denote the collection of 4 pitch classes that are held invariant: $\vb{W} = [0,2,4,6]$ (4-21, or the \textbf{W}hole tone tetramirror), $\vb{A} = [0,2,4,8]$ (4-24, or the \textbf{A}ugmented seventh chord), and $\vb{F} = [0,2,6,8]$ (4-25, or the \textbf{F}rench sixth set). The letter that appears in the superscript (outside of the parentheses) of the $\vb{S}$-type transformation denotes the set of pitches within the mystic/Wozzeck chord that does not change when the transformation is applied. The number within the parentheses, as with the $n = 4$ case, denotes the interval that is shifted.
\begin{figure*}[t]
  \centering
  \includegraphics[width=0.8\textwidth]{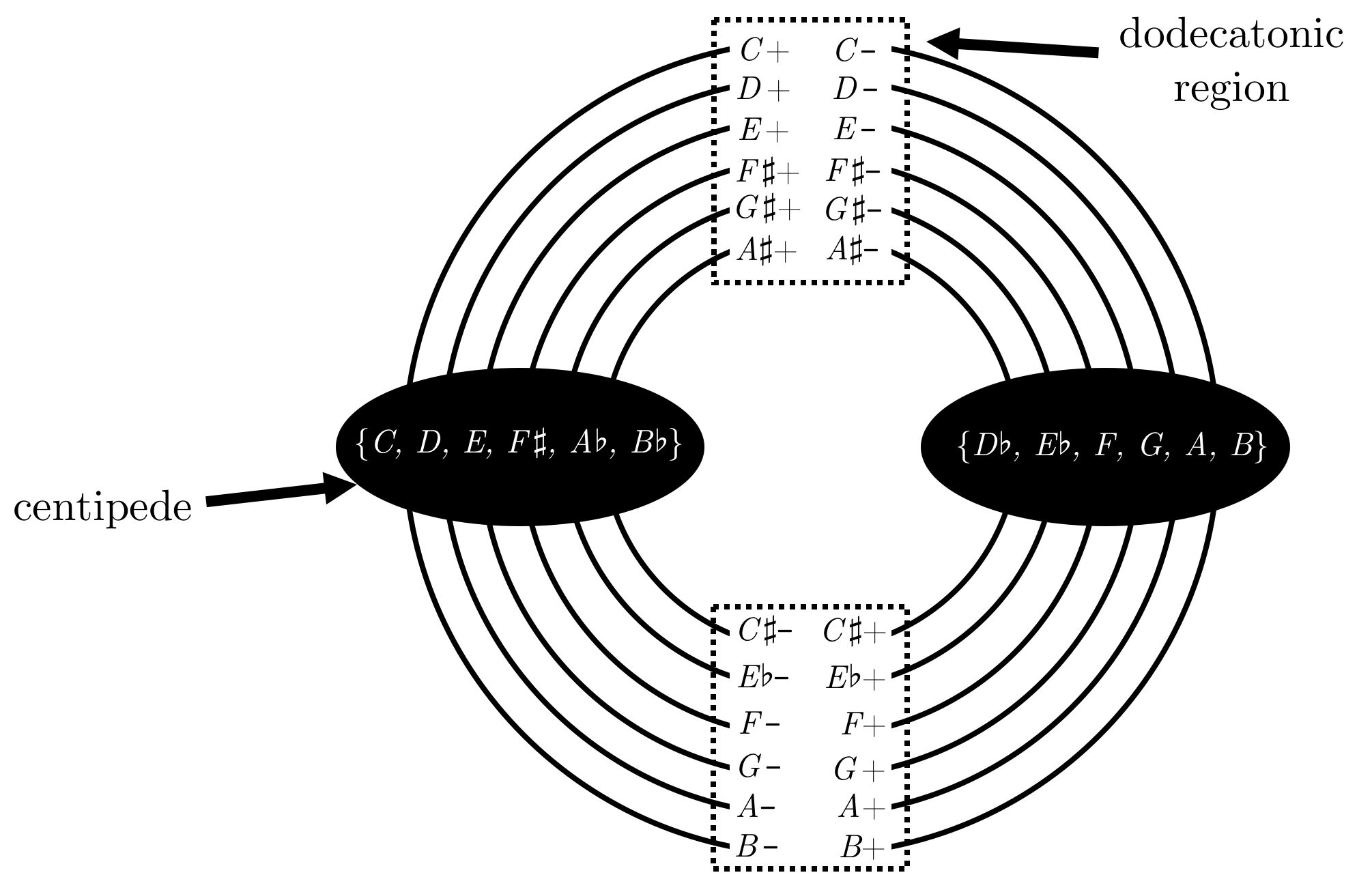}
  \caption{A unified voice-leading model for nearly symmetric hexachords.}
\end{figure*}

[5.8] As an example, suppose one applies $\vb{S^{A(3)}}$ to the Wozzeck chord $C+$, which is the pitch class collection $\{C, D\flat, E, F\sharp, G\sharp, B\flat\}$. The 2 augmented seventh chords which are subsets of this collection are $\{C, E, G\sharp, B\flat\}$ and $\{A\flat, C, E, G\flat\}$. Since the interval 3---a minor 3rd---is being shifted, the relevant augmented seventh chord that remains invariant during the transformation is $\{A\flat, C, E, G\flat\}$, and the minor 3rd of $B\flat$ and $D\flat$ is the interval that slides upward to $B$ and $D$. The new chord is formed by the pitch class collection $\{C, D, E, F\sharp, G\sharp, B\}$, which is $B-$, or the B mystic chord. The remaining Neo-Riemannian transformations in this region function the same way, and $\vb{S^{F}}$ is thereby an abbreviation of $\vb{S^{F(5)}}$.

[5.9] In \textbf{Figure 15}, I propose a unified voice-leading model analogous to Cohn's (2012) diagrams for connecting the Weitzmann waterbugs and Boretz spiders. By matching the root names of one centipede with the corresponding root names on the other centipede, 2 bridge regions arise between the centipedes. In analogy with hexatonic and octatonic regions, I refer to the boxes in \textbf{Figure 15} as dodecatonic regions. In \textbf{Figure 16}, I propose a visual representation of the voice-leading possibilities within the 2 dodecatonic regions, in analogy with \textbf{Figure 5} and \textbf{Figure 11}. Presumably, 5 spatial dimensions would be required to construct a convex polytope that presents all voice-leading possibilities within a dodecatonic region. Thus, \textbf{Figure 16} shows only a 2-dimensional reductive representation, constructed using the same logical procedure as I used in the previous section to ``flatten'' Childs' (1998) cubic network for the octatonic region.

[5.10] Any two chords that are connected by a solid line in the \textbf{Figure 16} are $P_{4,0}$-related. Thus, there are 5 available Neo-Riemannian transformations that can be used to transform a starting chord to the 5 $P_{4,0}$-related chords of the opposite modality. Using the naming convention introduced for the centipede's Neo-Riemannian transformations, these 5 ``slide'' transformations are: $\vb{S^1}$, $\vb{S^{3(A)}}$, $\vb{S^{3(W)}}$, $\vb{S^{5(A)}}$, and $\vb{S^{5(F)}}$. The number that appears first in the superscript denotes the interval between the 2 notes that are held invariant. The letter within the parentheses denotes the 4-note collection that is shifted. The letters follow the naming scheme described earlier in this section for the centipede transformations. For example, the transformation $\vb{S^1}$---which is an abbreviation of $\vb{S^{1(W)}}$---would transform $C+$ to $C-$, since the minor 2nd of pitch classes $C$ to $D\flat$ would be held invariant while the whole-tone tetramirror $\{E, F\sharp, A\flat, B\flat\}$ would be shifted downward a semitone.
\begin{figure*}[t]
  \centering
  \includegraphics[width=0.8\textwidth]{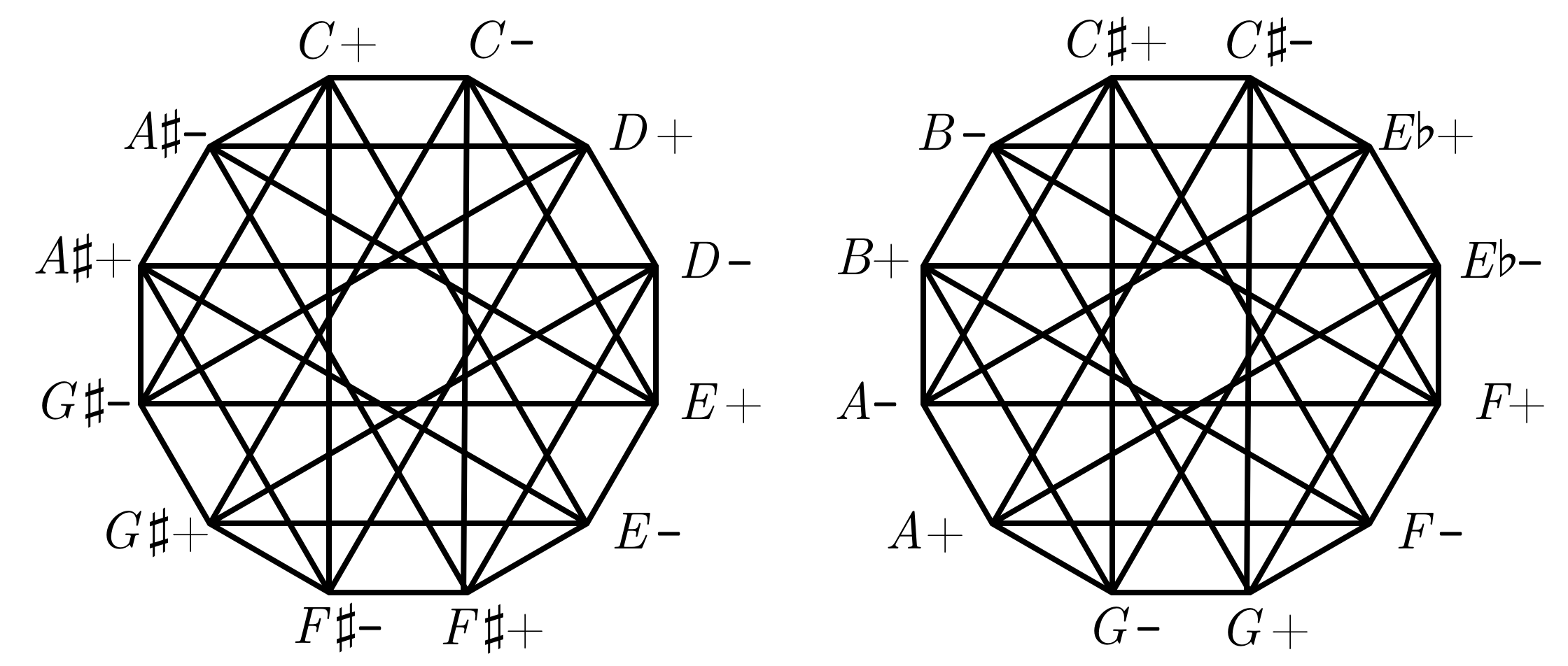}
  \caption{Dodecatonic regions, shown as 2-dimensional graphs.}
\end{figure*}
\begin{figure*}[t]
  \centering
  \includegraphics[width=0.4\textwidth]{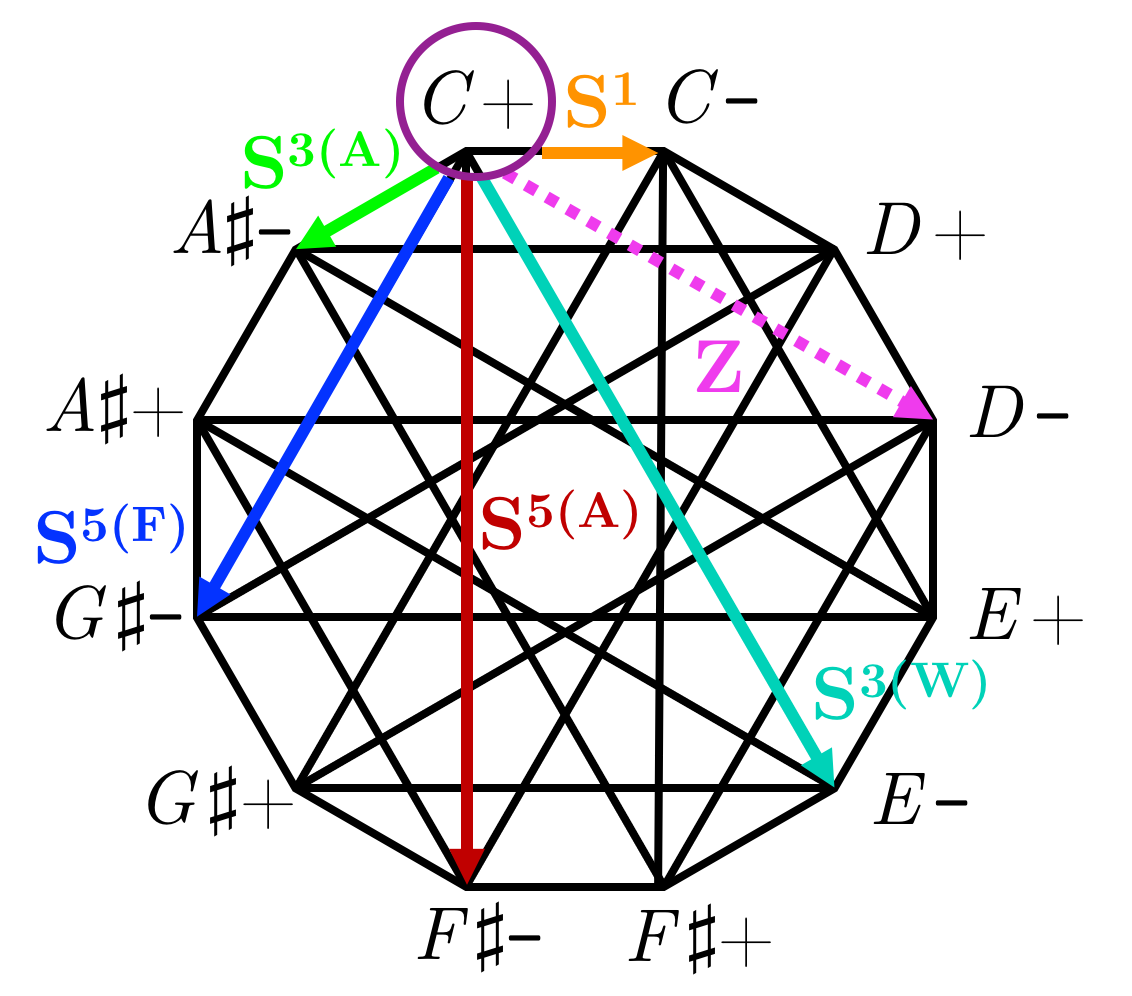}
  \caption{Available voice-leading transformations in a dodecatonic region.}
\end{figure*}

[5.11] Given a starting chord (suppose $C+$), there is one chord (in this case, $D-$) of the opposite modality that is not accessible via one of the 5 $\vb{S}$-type transformations; thus, it is not $P_{4,0}$-related to the starting chord. This chord shares no pitch classes with the starting chord, so functionally it is similar to the hexatonic pole for triads and the octatonic pole for seventh chords; thus, it can be referred to as the dodecatonic pole. The dodecatonic pole is accessible by the $\vb{Z}$ transformation\footnote[8]{Neo-Riemannian theorists often use $\vb{D}$ to denote motion to/from the ``dominant,'' so I use $\vb{Z}$, which stands for the German \textit{zw\"olf}, or ``twelve.''}. The geometric result of applying the $\vb{Z}$ transformation as well as the functions of the 5 $\vb{S}$-type transformations are shown in \textbf{Figure 17}.

[5.12] Within the dodecatonic region, one can construct various maximally smooth cycles. Sequences of $\vb{S}$-type transformations that generate closed paths with no loops along the solid lines in \textbf{Figure 16} are, indeed, maximally smooth cycles. The union of all pitch classes contained in the chords of any one of these maximally smooth cycles gives the full chromatic set 12-1, so these maximally smooth cycles can be called dodecatonic cycles.
 
\section*{VI. Discussion and Summary}
[6.1] In the previous 3 sections, I have shown that voice-leading models for nearly symmetric chords of cardinalities $n = 3$, $n = 4$, and $n = 6$ follow a set of patterns. In deriving the available Neo-Riemannian transformations for each type of chord, a symmetric chord is perturbed both downward and upward. Thus, the two chords that are generated from the perturbation of the same note are inversionally related.

[6.2] As I have mentioned before, the $(+)$ and $(-)$ labels assigned to major and minor, dominant seventh and half-diminished seventh, and Wozzeck and mystic chords signify the same direction of perturbation with regard to the symmetric chord from which each collection is generated. All 24 $(+)$ and $(-)$ chords can be organized into different groups, and one sees that the ``arthropod'' regions as well as the ``bridge'' regions in the unified voice-leading models are simply different ways of grouping these 24 chords in order to optimize voice-leading.

[6.3] If one chooses a particular chord from Weitzmann waterbug, Boretz spider, or centipede, one will always be able to find exactly 1 $P_{0,1}$-related chord of the opposite modality, and $n - 1$ $P_{2,0}$-related chords of the opposite modality within the same arthropod (Cohn 2012). If one chooses a particular chord from the hexatonic region, octatonic region, or dodecatonic region, one will always be able to find exactly $n - 1$ $P_{n-2,0}$-related chords of the opposite modality within the same region, and there will be exactly 1 ``polar'' chord which is the complement of the starting chord with respect to the union of all pitch classes within the region (Cohn 1996; Cohn 2012).
\begin{table*}[t]
\small
\caption{Summary of transformations, chord types, and voice-leading regions for nearly symmetric chords; adaptation and extension of Cohn (2012). Some row names from Cohn (2012) have been ommitted or modified, and I have added the $n = 6$ column.}
\centering
\begin{tabular}{|c|l|c|c|c|}
\hline
& Genus & Species, $n = 3$ & Species, $n = 4$ & Species, $n = 6$\\
\hhline{|=|=|=|=|=|}
1 & Symmetric partition & Augmented triad & Fully-diminished seventh & Whole-tone scale \\ \hline
2 & Downward SSD & Major triad & Dominant seventh & Wozzeck chord \\ \hline
3 & Upward SSD & Minor triad & Half-diminished seventh & Mystic chord \\ \hline
4 & Union of (2) and (3) & Consonant triads & Tristan genus & Mystic-Wozzeck genus \\ \hline
5 & arthropod region & Weitzmann waterbug & Boretz spider & Centipede (this paper) \\ \hline
6 & Bridge regions between (5)'s & Hexatonic region & Octatonic region & Dodecatonic region \\ \hline
7 & Transformations within (5) & \multicolumn{3}{|c|}{$\sim\sim\sim\sim\sim\sim\sim\sim\sim\sim\sim\sim\sim$} \\ \hline
& \;\textbullet\;For $P_{0,1}$-related chords & $\vb{R}$ & $\vb{R^*}$ & $\vb{R^{**}}$ \\ \hline
& \;\textbullet\;For $P_{2,0}$-related chords & $\vb{S}$, $\vb{N}$ & $\vb{S^{3(4)}}$, $\vb{S^{3(2)}}$, $\vb{S^{6}}$ & $\vb{S^{A(3)}}$, $\vb{S^{A(5)}}$, $\vb{S^{F}}$, $\vb{S^{W(1)}}$, $\vb{S^{W(3)}}$ \\ \hline
8 & Transformations within (6) & \multicolumn{3}{|c|}{$\sim\sim\sim\sim\sim\sim\sim\sim\sim\sim\sim\sim\sim$} \\ \hline
& \;\textbullet\;For $P_{n-2,0}$-related chords & $\vb{P}$, $\vb{L}$ & $\vb{S^2}$, $\vb{S^4}$, and $\vb{S^5}$ & $\vb{S^1}$, $\vb{S^{3(A)}}$, $\vb{S^{3(W)}}$, $\vb{S^{5(A)}}$, $\vb{S^{5(F)}}$ \\ \hline
& \;\textbullet\;For polar relation & $\vb{H}$ & $\vb{O}$ & $\vb{Z}$ \\ \hline
9 & Union of pitches in (6) & 6-20 & 8-28 & 12-1 \\ \hline
\end{tabular}
\end{table*}

[6.4] The union of all pitch classes within a particular hexatonic, octatonic, or dodecatonic region is also of set-theoretic interest. As Cohn (1996) mentions, within a given hexatonic region, there are only 6 unique pitch classes which are found; this collection is one of 4 distinct pitch-class sets generated from the set class 6-20. Each of the 4 pitch-class sets corresponds to a hexatonic region. Similar analysis shows that the octatonic regions correspond to set class 8-28 (Douthett and Steinbach, 1998). The dodecatonic regions correspond to the full chromatic set, 12-1.

[6.5] It is also of interest to note that the Neo-Riemannian transformations that are available in the arthropod regions are complementary to the transformations available in the bridge regions: the collections that are held invariant in one type of region are shifted in the other. For example, for the $n = 3$ case, the $\vb{S}$ holds a single note invariant while shifting the perfect 5th; the $\vb{P}$ transformation holds the perfect 5th invariant while shifting a single note. This relationship is easier to see in $n = 4$ and $n = 6$, as the names suggest a complementary relationship. As an example, one sees that for $n = 4$, the $\vb{S^{3(4)}}$ transformation acts within the arthropod region while $\vb{S^4} = \vb{S^{4(3)}}$ acts in the bridge region (Childs 1998). In general, one of the arthropod region transformations that connects $P_{2,0}$-related chords will have a complement in the bridge region, where transformations connects $P_{n-2,0}$-related chords.

[6.6] The Neo-Riemannian transformations, voice-leading regions, and set-theoretic properties I have discussed throughout this paper are summarized in \textbf{Table 1}, which is an extended version of Cohn's (2012) table covering $n = 3$ and $n = 4$. A further extension of the theory for $n = 6$ presented in this paper includes the rigorous development of a 5-dimensional \textit{Tonnetz} for voice-leading between any mystic and Wozzeck chords.

\section*{References}
{\small
\begin{hangparas}{1cm}{1}
\hangindent=1cm Childs, Adrian. 1998. ``Moving beyond Neo-Riemannian Triads: Exploring a Transformational Model for Seventh Chords.'' \textit{Journal of Music Theory} 42, no. 2: 181-193.

Cohn, Richard. 1996. ``Maximally Smooth Cycles, Hexatonic Systems, and the Analysis of Late-Romantic Triadic Progressions.'' \textit{Music Analysis} 15, no. 1: 9-40.

Cohn, Richard. 2000. ``Weitzmann's Regions, My Cycles, and Douthett's Dancing Cubes.'' \textit{Music Theory Spectrum} 22, no. 2: 89-103.

Cohn, Richard. 2012. \textit{Audacious Euphony: Chromatic Harmony and the Triad's Second Nature}. 2nd Edition. New York: Oxford University Press.

Douthett, Jack and Peter Steinbach. 1998. ``Parsimonious Graphs: A Study in Parsimony, Contextual Transformations, and Modes of Limited Transposition.'' \textit{Journal of Music Theory} 42, no. 2: 241-263.

Michelbacher, Abraham E. 1938. ``The biology of the garden centipede, \textit{Scutigerella immaculata}.'' \textit{Hilgardia} 11, no. 3: 55-148.

Tymoczko, Dmitri. 2011. \textit{A Geometry of Music}. New York: Oxford University Press.
\end{hangparas}
}
\end{document}